\documentclass[12pt,a4paper]{article}
\usepackage{epsfig,latexsym,amsfonts,amssymb,amsmath,amscd,
epic,graphics,theorem}
\theoremstyle{plain}
\newtheorem{lemma}{Lemma}

\newtheorem{theorem}[lemma]{Theorem}

{\theorembodyfont{\rmfamily}  \font\ncsc=cmcsc10
 \font\ntt=cmtt12

\begin{document}
\newcommand{\pperp}{\hbox{$\perp\hskip-6pt\perp$}}
\newcommand{\ssim}{\hbox{$\hskip-2pt\sim$}}
\newcommand{\aleq}{{\ \stackrel{3}{\le}\ }}
\newcommand{\ageq}{{\ \stackrel{3}{\ge}\ }}
\newcommand{\aeq}{{\ \stackrel{3}{=}\ }}
\newcommand{\bleq}{{\ \stackrel{n}{\le}\ }}
\newcommand{\bgeq}{{\ \stackrel{n}{\ge}\ }}
\newcommand{\beq}{{\ \stackrel{n}{=}\ }}
\newcommand{\cleq}{{\ \stackrel{2}{\le}\ }}
\newcommand{\cgeq}{{\ \stackrel{2}{\ge}\ }}
\newcommand{\ceq}{{\ \stackrel{2}{=}\ }}
\newcommand{\N}{{\mathbb N}}
\newcommand{\A}{{\mathbb A}}
\newcommand{\K}{{\mathbb K}}
\newcommand{\Z}{{\mathbb Z}}
\newcommand{\R}{{\mathbb R}}
\newcommand{\C}{{\mathbb C}}
\newcommand{\Q}{{\mathbb Q}}
\newcommand{\PP}{{\mathbb P}}
\newcommand{\mnote}{\marginpar}
\newcommand{\Id}{{\operatorname{Id}}}\newcommand{\Sym}{{\operatorname{Sym}}}
\newcommand{\oeps}{{\overline\eps}}
\newcommand{\oDel}{{\widetilde\Del}}
\newcommand{\real}{{\operatorname{Re}}}
\newcommand{\conv}{{\operatorname{conv}}}
\newcommand{\Span}{{\operatorname{Span}}}
\newcommand{\Ker}{{\operatorname{Ker}}}
\newcommand{\Hyp}{{\operatorname{Hyp}}}
\newcommand{\Fix}{{\operatorname{Fix}}}
\newcommand{\sign}{{\operatorname{sign}}}
\newcommand{\Tors}{{\operatorname{Tors}}}
\newcommand{\alg}{{\operatorname{alg}}}
\newcommand{\oi}{{\overline i}}
\newcommand{\oj}{{\overline j}}
\newcommand{\ob}{{\overline b}}
\newcommand{\os}{{\overline s}}
\newcommand{\oa}{{\overline a}}
\newcommand{\oy}{{\overline y}}
\newcommand{\ow}{{\overline w}}
\newcommand{\ot}{{\overline t}}
\newcommand{\oz}{{\overline z}}
\newcommand{\eps}{{\varepsilon}}
\newcommand{\proofend}{\hfill$\Box$\bigskip}
\newcommand{\Int}{{\operatorname{Int}}}
\newcommand{\pr}{{\operatorname{pr}}}
\newcommand{\grad}{{\operatorname{grad}}}
\newcommand{\rk}{{\operatorname{rk}}}
\newcommand{\im}{{\operatorname{Im}}}
\newcommand{\sk}{{\operatorname{sk}}}
\newcommand{\const}{{\operatorname{const}}}
\newcommand{\Sing}{{\operatorname{Sing}}\hskip0.06cm}
\newcommand{\conj}{{\operatorname{Conj}}}
\newcommand{\Cl}{{\operatorname{Cl}}}
\newcommand{\Crit}{{\operatorname{Crit}}}
\newcommand{\Ch}{{\operatorname{Ch}}}
\newcommand{\discr}{{\operatorname{discr}}}
\newcommand{\Tor}{{\operatorname{Tor}}}
\newcommand{\Conj}{{\operatorname{Conj}}}
\newcommand{\vol}{{\operatorname{vol}}}
\newcommand{\defect}{{\operatorname{def}}}
\newcommand{\codim}{{\operatorname{codim}}}
\newcommand{\tmu}{{\C\mu}}
\newcommand{\ov}{{\overline v}}
\newcommand{\ox}{{\overline{x}}}
\newcommand{\bw}{{\boldsymbol w}}
\newcommand{\bx}{{\boldsymbol x}}
\newcommand{\bu}{{\boldsymbol u}}
\newcommand{\bz}{{\boldsymbol z}}
\newcommand{\tet}{{\theta}}
\newcommand{\Del}{{\Delta}}
\newcommand{\bet}{{\beta}}
\newcommand{\kap}{{\kappa}}
\newcommand{\del}{{\delta}}
\newcommand{\sig}{{\sigma}}
\newcommand{\alp}{{\alpha}}
\newcommand{\Sig}{{\Sigma}}
\newcommand{\Gam}{{\Gamma}}
\newcommand{\gam}{{\gamma}}
\newcommand{\Lam}{{\Lambda}}
\newcommand{\lam}{{\lambda}}
\newcommand{\SC}{{SC}}
\newcommand{\MC}{{MC}}
\newcommand{\nek}{{,...,}}
\newcommand{\cim}{{c_{\mbox{\rm im}}}}
\newcommand{\clM}{\tilde{M}}
\newcommand{\clV}{\bar{V}}

\title{Appendix to ``Welschinger invariant and enumeration of\\
real rational curves"}
\author{Ilia Itenberg
\and Viatcheslav Kharlamov \and Eugenii Shustin}
\date{}
\maketitle

\begin{figure}[htb]
\hfill\includegraphics[height=1.2cm,width=10cm,angle=0,draft=false]{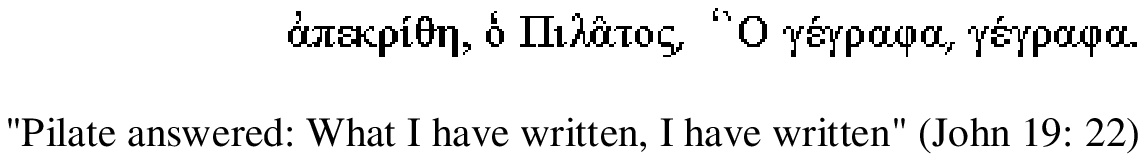}
\end{figure}

\section{Prologue}\label{introduction}

Since the appearance of our note \cite{IKS}, we received 
several
requests for a down-to-earth explanation of the invariance of
the Welschinger numbers, an explanation which would avoid 
symplectic geometry technology required by
a general symplectic setting
chosen in the original Welschinger's proof, see~\cite{W1}. Such an
explanation is certainly possible in the framework of the
classical algebraic geometry, especially, if one restricts itself
to the unnodal Del Pezzo surfaces. This class of surfaces includes
all the surfaces which are subject of \cite{IKS}.

In the proof given in this Appendix we follow the strategy of
Welschinger's proof. However, several technical details are
presented in another way because of our choice to work with
equations instead of parametrizations. Such an approach can be
useful for a study of higher genus cases in the algebraic geometry
framework. For this reason, we give preferences to proofs which
could be extended to higher genus even if they are a bit longer
than proofs based on parametrizations.

An {\it unnodal Del Pezzo surface} is a complex surface whose
anticanonical divisor is ample. In what follows we fix such an
unnodal Del Pezzo surface $\Sig$ equipped with a real structure
(i.e., an antiholomorphic involution) $c:\Sig\to\Sig$ and an ample
real divisor $D\subset\Sig$ such that $-D K_\Sig - 2
> 0$.

As is known, the set $R^0=R^0(D)$ of irreducible rational curves
in the linear system $|D|$ is locally closed and its closure
$R=R(D)$ is a projective variety. They are both of pure dimension
$r =-DK_\Sig-1$\footnote{It seems that the only explicit (and not
confused) formulation of the non-emptiness of $R^0(D)$ (for a nef
$D$) with a complete proof is found in \cite{GLS}, though it may
be deduced from the general Gromov-Witten theory; the count of the
local dimensions is by nowadays a routine application of the
deformation theory, which can be found almost everywhere.} and are
smooth at any $C\in |D|$ representing a nodal rational curve. The
Brusotti-Severi theorem (see \cite{Bru21}
and~\cite{Se})\footnote{This result was used by L.~Brusotti for a
construction of independent real variations of real nodes; note
that Brusotti treated as well the case of reducible nodal curves;
for far going contemporary generalizations see, for instance,
\cite{GrK}, Theorem 6.1(ii).} implies the following statements:
the tangent space to $R(D)$ at $C\in R^0(D)$ is the space
$T_CR=\Lambda(\Sing C)$ of curves in $|D|$ passing through $\Sing
C$, where $\Sing C$ is the set of singular points of~$C$; the
conditions imposed by different singular points are independent,
so $\codim\Lambda(\Sing C)$ is equal to the number $|\Sing C|$ of
points in $\Sing C$; and any linear system of dimension
$\codim\Lambda(\Sing C)$ intersecting $\Lambda(\Sing C)$ at $C$
transversely, induces a joint versal deformation of all the
singular points of $C$. Note also that according to the adjunction
formula one has
\begin{equation}
D^2 = r + 2|\Sing C| - 1\ . \label{e-new}\end{equation}

Let us fix an integer $m$ such that $0\le 2m\le r$ and introduce a
real structure $c_{r,m}$ on $\Sig^r$ which maps $(z_1, \ldots ,
z_r)\in\Sig^r$ to $(z'_1, \ldots , z'_r)\in\Sig^r$ with $z'_i =
c(z_i)$ if $i>2m$, and $(z'_{2j-1}, z'_{2j}) = (c(z_{2j}),
c(z_{2j-1}))$ if $j\le m$. With respect to this real structure, a
point $\bw=(z_1,...,z_r)$ is real, {\it i.e.},
$c_{r,m}$-invariant, if and only if $z_i$ belongs to the real part
$\R\Sig$ of $\Sig$ for $i>2m$ and $z_{2j-1},z_{2j}$ are conjugate
for $j \le m$. In what follows we work with an open dense subset
$\Omega_{r,m}(\Sig)$ of $\R\Sig^r=\Fix\, c_{r,m}$ constituted of
$c_{r,m}$-invariant $r$-tuples $\bw=(z_1,...,z_r)$ with pairwise
distinct $z_i\in\Sig$.

By abuse of language, we say that a curve $C$ in $\Sig$ passes
through $\bw\in\Sig^r$ if $C$ contains all the components
$z_i\in\Sig$ of $\bw$. If $C$ is nonsingular at each of the
components $z_i$ of $\bw$, we say that $C$ is nonsingular in
$\bw$. Moreover, when it can not lead to a confusion we denote by
the same symbol $\bw$ the set of its components.

In the spirit of \cite{W,W1}, given a {\bf generic}
$\bw\in\Omega_{r,m}(\Sig)$, we count the number
$N_{r,m}^{even}(\bw)$ (resp., $N_{r,m}^{odd}(\bw)$) of irreducible
real rational curves in $|D|$ passing through $\bw$ and having
even (resp., odd) number of solitary nodes. (Note that for a
generic $\bw$ every rational curve in $|D|$ passing through $\bw$
is irreducible and nodal.) The {\it Welschinger number}
$W_{r,m}(\bw)$ is defined by
$W_{r,m}(\bw)=N_{r,m}^{even}(\bw)-N_{r,m}^{odd}(\bw)$.\footnote{In
\cite{IKS} we considered only the case $m=0$ and used a slightly
different notation.}

\begin{theorem}\label{p1}
{\rm (J.-Y.~Welschinger, see~\cite{W, W1})}. The value
$W_{r,m}(\bw)$ does not depend on the choice of a {\rm
(}generic{\rm )} element $\bw$ in a given connected component of
$\Omega_{r,m}(\Sig)$. In particular, if $\R\Sig$ is connected, and
$\bw_1$ and $\bw_2$ are two generic elements of
$\Omega_{r,m}(\Sig)$, then $W_{r,m}(\bw_1) = W_{r,m}(\bw_2)$.
\end{theorem}

In a connected component of $\Omega_{r,m}(\Sig)$, we can join any
two generic elements by a sequence of paths such that
\begin{itemize}
\item either one of the real
components of $\bw$, say $z_r$, travels in $\R\Sig$ along a real
segment of a smooth generic real algebraic curve avoiding
$z_1,...,z_{r-1}$ and crossing all $(-1)$-curves transversely,
whereas the points $z_1,...,z_{r-1}$ stay fixed,
\item or a pair of imaginary conjugate
components of $\bw$, say, $z_1,z_2$, travel along generic smooth
conjugate arcs in the non-real part $\Sig$, avoiding $z_3,...,z_r$
and all $(-1)$-curves, whereas the components $z_3,...,z_r$ stay
fixed.
\end{itemize}

We intend to show that $W_{r,m}(\bw)$ remains constant along any
of such paths.

\section{Moving a real point of the configuration}\label{proof}

We will use the following observation.

\vskip5pt

\noindent
{\bf Transfer of genericity.} {\it Let $Y$ be an irreducible
complex variety equipped with a real structure. Assume that the
real part $\R Y$ of $Y$ contains nonsingular points. Then any
generic point of $\R Y$ is also a generic point of $Y$}.

\vskip5pt

In particular, a generic point of $\Omega_{r,m}(\Sig)$ is a
generic point of $\Sig^r$.

Assume that $2m<r$. Let
$\bw'=(z_1,...,z_{r-1})\in\Omega_{r-1,m}(\Sig)$ be generic, and
let the point $z_r$ move along a segment $\sig$ of a real part of
some generic real smooth algebraic curve $S_0$ in $\Sig$. In
particular, we suppose that $\sig$ does not contain any component
of $\bw'$. Denote by $V(\bw')\subset|D|$ {\it the set of
irreducible rational nodal curves in $|D|$ which pass through
$\bw'$ and are nonsingular at $\bw'$.} Denote by $\clV(\bw')$ the
closure of $V(\bw')$, by $\R\clV(\bw')$ the real part of
$\clV(\bw')$, and by $\R{V(\bw')}$ the real part of ${V(\bw')}$.

Let $\Lam(z_r)\subset|D|$ be the linear system of curves passing
through $z_r$. We study bifurcations of the curves in the set
$\Lam(z_r)\cap\R \clV(\bw')$ along a move of $z_r$ through $\sig$,
and show that $W_{r,m}(\bw)$ does not change in these
bifurcations.

\begin{lemma}\label{l1}
The set $V(\bw')$ {\rm (}respectively, $\R V(\bw')${\rm )} is a
smooth one-dimensional subvariety of $|D|$ {\rm (}respectively,
$\R |D|$ {\rm ).} Moreover, the intersection of $R$ with
$\Lambda(\bw')$ is transversal along $V(\bw')$, and in its turn
$V(\bw')$ at a point $C\in V(\bw')$ intersects transversely with a
linear system $\Lam(z_r)$, if $z_r\not\in\Sing C$.
\end{lemma}

{\bf Proof.} All the claims follow from the fact that, for any
$C\in V(\bw')$ and any $z_r\not\in\Sing C$, the points in $\bw
\cup \Sing C$ impose independent conditions on the curves of
$|D|$. The latter statement means that
$\Lambda(\bw)\cap\Lambda(\Sing C)=\{C\}$ and follows, for
instance, from the B\'ezout theorem, since $C$ is irreducible and
$r+2|\Sing C|=D^2+1>D^2$ (see~(\ref{e-new})).
\proofend

In particular, Lemma~\ref{l1} implies that $W_{r,m}(\bw)$ may
change only when $\Lam(z_r)$ crosses elements of
$\R\clV(\bw')\backslash \R V(\bw')$ or becomes tangent to $\R
V(\bw')$, and that such critical situations appear for at most
finitely many $z_r$ in $\sig$.

\medskip

\begin{lemma}\label{l2}
The elements of $\clV(\bw')\backslash V(\bw')$ are
\begin{enumerate}
\item[{\rm (i)}]
irreducible rational curves whose collection of singular points
consists of a cusp $A_2$ and $D(D+K_\Sig)/2$ nodes,
\item[{\rm (ii)}]
irreducible rational curves whose collection of singular points
consists of a tacnode $A_3$ and $D(D+K_\Sig)/2-1$ nodes,
\item[{\rm (iii)}]
irreducible rational curves whose collection of singular points
consists of a triple point $D_4$ and $D(D+K_\Sig)/2-2$ nodes,
\item[{\rm (iv)}]
reduced reducible curves, splitting into two rational nodal
curves, intersecting transversely and only at their nonsingular
points,
\item[{\rm (v)}]
irreducible rational nodal curves with $D(D+K_\Sig)/2 + 1$ nodes
such that exactly one of the nodes coincides with a point of
$\bw'$.
\end{enumerate}
Moreover, in the first four cases, the curves are nonsingular at
$\bw'$.
\end{lemma}

{\bf Proof.} According to the
observation
on the transfer of genericity, it
is sufficient to prove the statement of Lemma for generic $\bw'
\in \Sig^{r-1}$. Notice that the linear system $\Lam(\bw')$ does
not meet subvarieties of $|D|$ of dimension $<r-1=-DK_\Sig-2$.
Indeed, choosing $z_1$ outside a curve $C\in X \subset|D|$, we
obtain $\dim(\Lam(z_1)\cap X)=\dim X-1$ and then proceed by
induction.

Thus, to show that the elements of $\clV(\bw')\backslash V(\bw')$
can not have more complicated singularities than those pointed in
(i)--(iv) it is sufficient to carry out suitable dimension counts.

Assume that $C\in \clV(\bw')$ is reduced irreducible. Then it is
rational. The Zariski tangent space to the germ at $C$ of the
equisingular stratum in $|D|$ is the projectivization of
$H^0(\Sig,{\cal J}_{Z^{es}(C)}(D))$, where ${\cal J}_{Z^{es}(C)}$
is the ideal sheaf of the zero-dimensional scheme $Z^{es}(C)$
supported at $\Sing(C)$ and defined at $z\in\Sing(C)$ by the
equisingular ideal (see \cite{DH,Wa}). The scheme $Z^{es}(C)$
contains the conductor scheme $Z^{cond}(C)$ defined at any
$z\in\Sing(C)$ by the conductor ideal, and, furthermore,
$$\deg Z^{es}(C)-\deg
Z^{cond}(C)=\sum_{p\in\Sing(C)}(\tau'(p,C)-\del(p,C))\ ,$$ where
$\tau'(p,C)$ is the codimension of the equisingular ideal, and
$\del(p,C)$ is the $\del$-invariant. In our case $\deg
Z^{cond}(C)=D(D+K_\Sig)/2+1$, and hence (see, for example,
\cite{DH}), $\deg Z^{es}(C)=D(D+K_\Sig)/2+2$ only in cases
(i)--(iv). For other collections of singularities, $\deg
Z^{es}(C)\ge D(D+K_\Sig)/2+3$. To finish our dimension count it
remains to show that, in the latter case, $h^0(\Sig,{\cal
J}_{Z^{es}(C)}(D))\le-DK_\Sig-2$. For this purpose, pick any
zero-dimensional scheme $Z$ of degree $D(D+K_\Sig)/2+3$ between
$Z^{es}(C)$ and $Z^{cond}(C)$. Then
$$h^0(\Sig,{\cal J}_{Z^{es}(C)}(D))\le h^0(\Sig,{\cal
J}_Z(D))=\frac{D(D-K_\Sig)}{2}+1-\deg Z+h^1(\Sig,{\cal J}_Z(D))$$
$$ =-DK_\Sig-2+h^1(\Sig,{\cal J}_Z(D))=-DK_\Sig-2\ .$$ The relation
\begin{equation}h^1(\Sig,{\cal J}_Z(D))=0\
,\label{enn2}\end{equation} which we use here, follows from
$C\supset Z^{es}(C)\supset Z$ and from the exact sequence
$$0=H^1(\Sig,{\cal O}_{\Sig})\to H^1(\Sig,{\cal J}_Z(D))\to
H^1(C,{\cal J}_{Z/C}\otimes{\cal O}_{\Sig}(D))=0\ .$$ In turn, the
latter vanishing comes from the Riemann-Roch theorem, since (see
\cite{ACGH}) in view of $Z\supset Z^{cond}(C)$, the sheaf ${\cal
J}_{Z/C}\otimes{\cal O}_{\Sig}(D)$ on $C$ lifts up to a sheaf
${\cal O}_{C^\nu}(D))$ of the normalization $C^\nu$, and
$$\deg{\cal O}_{C^\nu}(D))=\deg{\cal J}_{Z/C}\otimes{\cal
O}_{\PP^2}(d)=D^2-\deg Z-\deg Z^{cond}(C)$$ $$\ge-DK_\Sig-4>-2=
2g(C)-2.$$

Assume now that $C\in\clV(\bw')\backslash V(\bw')$ splits into
irreducible components $C_1,...,C_k$, $k\ge 2$ (may be,
coinciding). The classical theorem (see \cite{N} for a modern
proof) states that $g(C_1)+...+g(C_k)\le g(C')=0$, where $C'$ is a
generic curve in $V(\bw')$. Hence, all $C_1,...,C_k$ are rational.
Since $-C_iK_\Sig -1\ge 0$ for any $i = 1, \ldots , k$ (because
the surface is unnodal), and since a rational curve $C_i$ cannot
contain more than $-C_iK_\Sig-1$ generic points, the curve $C$
does not contain more than $-DK_\Sig-k$ generic points. Thus, $k =
2$ and the rational curves $C_1$ and $C_2$ pass through
$-C_1K_\Sig-1$ and $-C_2K_\Sig-1$ generic points of $\Sig$,
respectively. Therefore, both the curves $C_1$ and $C_2$ are nodal
and intersect each other transversely.

The last statement of Lemma also follows from a dimension count.
Namely, we deduce from the preceding dimension count that the
curves of the types described in the first four cases and having a
singularity at a fixed generic point form a variety of dimension
$\le (-D K_\Sig - 2) - 2$. Since the number of additional points
to mark is equal to $-D K_\Sig - 3$, for generic marked points
such curves do not appear.

According to a similar count, the curves with nodes at two marked
points do not appear if the marked points are generic.\proofend

Let $\cal T$ be the curve generated by the double points of $C\in
V(\bw')$, and $\cal T'$ be the underlying reduced curve. To
control the behavior of Welschinger number at critical moments we
put special attention, in accordance with Lemmas \ref{l1} and
\ref{l2}, to the curves $C\in \clV(\bw')\backslash V(\bw')$ and
$\cal T$. Namely, we impose an {\bf extra condition} on $S_0$
requiring that it crosses $\cal T'$ and all $C\in
\clV(\bw')\backslash V(\bw')$ transversely and at their generic
nonsingular points.

\begin{lemma}\label{lT}
The number $W_{r,m}(\bw)$ does not change when $z_r$ crosses $\cal
T$ along $\sig$.
\end{lemma}

{\bf Proof.} Let $z_r^0$ be such a crossing point. It does not
belong to any $C\in \clV(\bw')\backslash V(\bw')$, as it follows
from our choice of $S_0$. Pick a curve $C^0\in \Lambda(z_r^0)\cap
V(\bw')$ with a node at $z_r^0$. The equality~(\ref{e-new})
implies that the nearby curves $C^t\in V(\bw')$ form a regular
homotopy covering which covers twice a neighborhood of $z_r^0$;
each of the two branches of $C^0$ at $z_r$ provides a leaf of this
trivial $2$-covering. The curve $C^0$ is real and the real parts
$\R C^t$ of the nearby $C^t\in \R V(\bw')$ form a regular
equivariant homotopy which covers twice a real neighborhood of
$z_r^0$ if $z_r^0$ is not a solitary point of $\R C^0$, otherwise
they cover exclusively $\R\cal T$. Hence, the real solutions with
$z_r$ on one side of $\R \cal T$ are real regular homotopy
equivalent to the respective real solutions with $z_r$ on the
other side of $\R \cal T$, and, therefore, $W_{r,m}(\bw)$ does not
change.
\proofend

\begin{lemma}\label{l3}
Let $C\in\R \clV(\bw')$ be as in the cases {\rm (i), (ii)}, or
{\rm (iii)} of Lemma \ref{l2}. Then $W_{r,m}(\bw)$ does not change
when $z_r$ crosses $C$ along $\sig$.
\end{lemma}

{\bf Proof.} Let us check, first, that the linear system
$\Lam(\bw')$ induces a joint versal deformation of all the
singular points of $C$. Such a property of $\Lam(\bw')$ is
equivalent to the transversality of intersection of $\Lam(\bw')$
and $|{\cal J}_{Z^{es}(C)}(D)|$ in $|D|$. As we have seen in the
proof of Lemma 3, $\dim|D|=D(D-K_\Sig)/2$,
$\dim\Lam(\bw')=\dim|D|-r+1$, and
$$\dim|{\cal J}_{Z^{es}(C)}(D)|=\dim|D|-\deg Z^{es}+h^1(\Sig,{\cal J}_{Z^{es}(C)}(D))$$
$$=\frac{D(D-K_\Sig)}{2}-\left(\frac{D(D+K_\Sig)}{2}+2\right)+
h^1(\Sig,{\cal
J}_{Z^{es}(C)}(D))\stackrel{\text{(\ref{enn2})}}{=}r-1\ .$$ Hence
it is sufficient to show that
\begin{equation}\dim\left(\Lam(\bw')\cap|{\cal
J}_{Z^{es}(C)}(D)|\right)=0\ \label{e5}\end{equation} under a
suitable choice of $\bw'$, and it remains to notice that choosing
one-by-one generic points $z_1,...,z_{r-1}$ on $C$ outside $C_1\in
|{\cal J}_{Z^{es}(C)}(D)|, C_2\in |{\cal
J}_{Z^{es}(C)}(D)|\cap\Lambda(z_1)$,\dots , we reduce the
dimension of $\Lam(\bw')\cap|{\cal J}_{Z^{es}(C)}(D)|$ up to zero.

On the other hand, as it follows from Lemma~\ref{l2}, the set
$\Sing(C)$ contains a singular point $z^*$ of type $A_2$, $A_3$ or
$D_4$, and $\Sing(C)\backslash\{z^*\}$ consists of
$D(D+K_\Sig)/2+2-\mu(z^*,C)$ nodes. For any node $q$ of $C$, the
germ $N_q(C)$ at $C$ of the set of curves in $|D|$ having a node
in a neighborhood of $q$, is a smooth hypersurface in $|D|$.
Therefore, the above transversality implies that the
$\mu(z^*,C)$-dimensional germ
$P=\Lam(\bw')\cap\bigcap_{q\in\Sing(C)\backslash\{z^*\}}N_q(C)$ at
$C$ is smooth, and induces a versal deformation of the singular
point $(C,z^*)$.

The germ of $\clV(\bw')$ at $C$ is contained in $P$ and coincides
with the $\del$-constant stratum in the corresponding versal
deformation. It is well-known that, for $z^*$ of type $A_2$, this
stratum has an ordinary cusp at $C$ (it is described by the
discriminant equation $\lam^3/27+\mu^2/4=0$ in the standard versal
deformation $x^3+\lam x+\mu$ of the cusp). For $z^*$ of type $A_3$
or $D_4$, the germ of $\clV(\bw')$ at $C$ is smooth by \cite{DH},
Proposition 4.17(2). Denote by $T$ the (one-dimensional) tangent
line to $\clV(\bw')$ at $C$. Let $z_r^0$ be an intersection point
of $\sig$ with $C$. The intersection of $P$ and $\Lam(z_r^0)$ is
transversal, and $P\cap\Lam(z_r^0)$ is a smooth line (resp,
surface, three-fold) in the surface (resp., three-fold, four-fold)
$P$, if $z^*\in C$ is of type $A_2$ (resp., $A_3$, $D_4$).
Furthermore, $T_CP\cap\Lam(z_r^0)$ intersects transversally with
the line $T$ in the tangent space $T_CP$ to $P$ at $C$. Indeed, in
view of the generic choice of $z_r^0$ on $C$, the latter
transversality simply means that $T$ contains a curve different
from $C$. Then, in particular, any linear system $\Lam(z_r)$ with
$z_r$ belonging to the germ of $\sig$ at $z_r^0$, crosses
$\clV(\bw')$ transversely in a neighborhood of $C$. Hence, varying
$z_r$ along the germ of $\sig$ at $z_r^0$, we observe for $z_r\ne
z_r^0$, that
\begin{itemize}
\item
if $z^*\in C$ is of type $A_2$, then $\R V(\bw')\cap\Lam(z_r)\cap
P$ is empty or consists of two points, one corresponding to a
curve with a solitary node in a neighborhood $U^*$ of $z^*$ in
$\R\Sig$, the other corresponding to a curve with a non-solitary
node in $U^*$;
\item
if $z^*\in C$ is a non-solitary tacnode $A_3$, then $\R
V(\bw')\cap\Lam(z_r)\cap P$ consists of one point, corresponding
either to a curve with two non-solitary nodes in $U^*$, or to a
curve without nodes in $U^*$;
\item
if $z^*\in C$ is a solitary tacnode $A_3$, then $\R
V(\bw')\cap\Lam(z_r)\cap P$ consists of one point, corresponding
either to a curve with two solitary nodes in $U^*$, or to a curve
without nodes in $U^*$;
\item
if $z^*\in C$ is of type $D_4$ with three real local branches,
then $\R V(\bw')\cap\Lam(z_r)\cap P$ consists of one point,
corresponding to a curve with three non-solitary nodes in $U^*$;
\item
if $z^*\in C$ is of type $D_4$ with one real local branch, then
$\R V(\bw')\cap\Lam(z_r)\cap P$ consists of one point,
corresponding to a curve with one solitary node in $U^*$.
\end{itemize}
In all the above cases, the Welschinger number $W_{r,m}(\bw)$ does
not change when $z_r$ crosses $C$.\proofend

\begin{lemma}\label{l4}
Let $C\in \R\clV(\bw')\backslash \R V(\bw')$ be as in case {\rm
(iv)} of Lemma \ref{l2}. Then $W_{r,m}(\bw)$ does not change when
$z_r$ crosses $C$ along $\sig$.
\end{lemma}

{\bf Proof.} According to Lemma~\ref{l2}, the curve $C$ splits
into two irreducible rational nodal curves $C_1$ and $C_2$ which
intersect each other transversely and only at their nonsingular
points. These curves are either both real or complex conjugate to
each other. Since by the Riemann-Roch theorem
$\dim\Lam(\bw'\cup\Sing(C))=0$, the germ of $ \R\clV(\bw')$ at $C$
is the union of $s$ transverse smooth branches, where $s$ is the
number of real intersection points of $C_1$ and $C_2$, the tangent
lines to these branches are
$T_q=\Lam(\bw'\cup\Sing(C)\backslash\{q\})$, $q\in C_1\cap
C_2\cap\R\Sig$, and when $C$ moves along such a branch all the
nodes are preserved, except $q$ which undergoes a Morse
transformation. By our convention on the choice of $\sig$, any
intersection point $z_r^0$ of $\sig$ and $C$ lies outside some
preselected generic curves in $T_q$, $q\in C_1\cap C_2
\cap\R\Sig$. Then, inside $\Lam(\bw')$, the linear system
$\Lam(\bw',z_r^0)$ intersects transversely the branches of
$\R\clV(\bw')$ at $C$. The number of solitary nodes (as well as
the number of non-solitary nodes) is not changing under a Morse
transformation outside the nodes; hence, the Welschinger number
$W_{ r,m}(\bw)$ does not change in this situation.\proofend

\begin{lemma}\label{l5}
Let $C\in \R\clV(\bw')\backslash \R V(\bw')$ be as in the case
{\rm (v)} of Lemma \ref{l3}. Then $W_{r,m}(\bw)$ does not change
when $z_r$ crosses $C$ along $\sig$.
\end{lemma}

{\bf Proof.} According to Lemma~\ref{l2}, we can assume that
$2m\le r-2$, and $C$ is an irreducible rational nodal curve with a
node at $z_{r-1}\in\R\Sig$ and nonsingular at the other points of
$\bw'$. The germ of the smooth surface
$P=\Lam(\bw')\cap\bigcap_{q\in\Sing(C)\backslash\{z_{r-1}\}}
N_q(C)$ (recall that $N_q(C)$ is the germ at $C$ of the set of
curves in $|D|$ having a node in a neighborhood of $q$) represents
a hyperplane section of a versal deformation of the singularity
$z_{r-1}\in C$, namely the section determined by the condition to
pass through $z_{r-1}$. The germ of $\clV(\bw')$ at $C$ lies in
$P$ and consists of two smooth branches, $B_1$ and $B_2$, each
representing an equisingular deformation. For each of the branches
pick a curve representing a point different of $C$ on the tangent
to the branch. By the genericity assumption, we may suppose that
$\sig$ crosses $C$ at a point $z_r^0$, which does not belong to
the above curves. Then the linear system $\Lam(\bw',z_r^0)$
crosses transversely the tangents to $B_1$ and $B_2$.

If $z_{r-1}$ is a solitary node of $\R C$, then $B_1$ and $B_2$
are imaginary, and hence $\Lam(\bw)\cap \R V(\bw')\cap
P=\emptyset$ as $z_r\ne z_r^0$. If $z_{r-1}$ is a non-solitary
node of $\R C$, then $B_1$ and $B_2$ are real, and hence
$\Lam(\bw)\cap \R V(\bw')\cap P$ consists of two elements, and
these elements correspond to real curves which are real
equisingular deformation equivalent to $C$.

Thus, $W_{r,m}(\bw)$ does not change.
\proofend

In the case $m = 0$, Theorem~\ref{p1} immediately follows from the
results of this Section. For $m>0$, it remains to prove the
invariance of $W_{r,m}$ under moves of a pair of conjugated marked
points.

\section{Moving a pair of imaginary points of the configuration.
End of the proof of Theorem~\ref{p1}}\label{proof2}

Assume that $m>0$.

Let $\bw'=\{z_3,...,z_r\}$ be a generic element in
$\Omega_{r-2,m-1}(\Sig)$, the point $z_1$ move along a generic
smooth path $\sig$ in $\Sig \backslash \R \Sig$, and $z_2=\bar
z_1$. Denote by $V(\bw')$ the set of irreducible rational nodal
curves in $|D|$ which pass through $\bw'$ and are nonsingular at
$\bw'$. Denote by $\clV(\bw')$ the closure of $V(\bw')$, by
$\R\clV(\bw')$ the real part of $\clV(\bw')$, and by $\R{V(\bw')}$
the real part of ${V(\bw')}$.

Let $\Lam(z_1,z_2)\subset|D|$ be the linear system of curves
passing through $z_1$ and $z_2$. We study bifurcations of the set
$\Lam(z_1,z_2)\cap\R\clV(\bw')$ along the path $\sig$.

\begin{lemma}\label{l1X}
The set $\R V(\bw')$ is a smooth two-dimensional subvariety of $\R
|D|$. Moreover, the intersection of $R$ with $\Lambda(\bw')$ is
transversal along $\R V(\bw')$, and in its turn $\R V(\bw')$
intersects transversely with the real part of any linear system
$\Lam(z_1,z_2)$, where $z_1\in\sig$, $z_2=\bar z_1$.
\end{lemma}

{\bf Proof.} The proof almost coincides with the proof of Lemma
\ref{l1}. The only new point is to show that, for generic $\bw'$
and $\sig$, the curves $C\in \R V(\bw')$ are nonsingular at $z_1$
and $z_2$. The latter statement follows from the fact that $\sig$
can be chosen in such a way that it avoids the singular points of
curves $C \in \R V(\bw')$, which is a generic condition on $\sig$.
\proofend

Lemma~\ref{l1X} implies that for all but finitely many positions
of $z_1,z_2$ in a path, one has
$\Lam(z_1,z_2)\cap\clV(\bw')\subset\Lam(z_1,z_2)\cap V(\bw')$, and
that $W_{r,m}(\bw)$ may change only when $\Lam(z_1,z_2)$ crosses
elements of $\R\clV(\bw')\backslash \R V(\bw')$.

\begin{lemma}\label{l2X}
When $z_2=\bar z_1$ and $z_1$ moves along $\sig$, the linear
system $\Lam(z_1,z_2)$ crosses only one-dimensional strata of
$\R\clV(\bw')\backslash \R V(\bw')$, whose elements are as
follows:
\begin{enumerate}
\item[(i)]
irreducible rational curves with a collection of singular points,
consisting of a cusp $A_2$ and $D(D+K_\Sig)/2$ nodes,
\item[(ii)]
irreducible rational curves with a collection of singular points,
consisting of a tacnode $A_3$ and $D(D+K_\Sig)/2-1$ nodes,
\item[(iii)]
irreducible rational curves with a collection of singular points,
consisting of a triple point $D_4$ and $D(D+K_\Sig)/2-2$ nodes,
\item[(iv)] reduced, reducible
curves, splitting into two rational nodal curves, intersecting
transversally and only at their nonsingular points,
\item[(v)]
irreducible rational nodal curves with $D(D+K_\Sig)/2 + 1$ nodes
such that exactly one of these nodes coincides with a real point
of $\bw'$.
\end{enumerate}
Moreover, in the first four cases the curves are nonsingular at
$\bw'$.
\end{lemma}

{\bf Proof.} By Lemma \ref{l2} the curves as in (i)--(v) form an
open dense set in $\clV(\bw')\backslash V(\bw')$. Denote it by
$Z$. The result follows as soon as we choose $\sig$ to avoid the
singular points of the curves $C\in Z$ as well as the points of
the curves in the finite set $(\clV(\bw')\backslash
V(\bw'))\backslash Z$.
\proofend

\medskip

Now we specify a particular genericity condition and {\bf require}
that the path $\sig$ crosses the curves, mentioned in Lemma
\ref{l2X}, only at their nonsingular points.

\begin{lemma}\label{l3X}
The value of $W_{r,m}(\bw)$ does not change when $z_1$ and $z_2$
cross a curve of any of the types described in Lemma \ref{l2X}
along the paths $\sig$ and $c(\sig)$.
\end{lemma}

{\bf Proof.} We use the transversality arguments as in the proofs
of Lemmas \ref{l3}, \ref{l4}, and \ref{l5}. The transversality is
assured by the above requirement.

Assume, for instance, that $\sig$ crosses a curve $C$ as described
in case (i), having a cusp at
$z^*\in\R\Sig\backslash\{z_1,...,z_r\}$. Since $z_1,\dots, z_r$
are nonsingular points of $C$ (the above requirement), the
three-dimensional smooth germ of the variety
$P=\Lam(\bw')\cap\bigcap_{q\in\Sing(C)\backslash\{z^*\}}N_q(C)$
represents a versal deformation of the singular point $z^*$ of
$C$, in which the discriminant is a surface with a cuspidal edge.
The cuspidal edge divides the real part of the discriminant into a
smooth surface corresponding to curves with a solitary node in a
neighborhood of $z^*$, and a smooth surface corresponding to
curves with a non-solitary node in a neighborhood of $z^*$. The
(projective) tangent plane to the discriminant at $C$ is naturally
isomorphic to the linear system $\Lam(\bw'\cup\Sing(C))$, and, by
the Riemann-Roch theorem, it intersects transversely the linear
system $\Lam(z^0_1,z^0_2)$, $z^0_1,z^0_2$ standing for the
travelling points $z_1,z_2$, corresponding to $C$. Hence, varying
$z_1,z_2$ along the germs of $\sig, c(\sig)$ at $z^0_1,z^0_2$,
respectively, we obtain a transversal intersection of the
discriminant with $\Lam(z_1,z_2)$, which then consists of two
imaginary curves, or of two real curves with the opposite
contributions to $W_{r,m}(\bw)$.

Similarly one proves the invariance of $W_{r,m}(\bw)$ in other
bifurcations. \proofend

Thus the proof of Theorem~\ref{p1} is completed.\proofend

\section{Epilogue}\label{Epilogue}
{\bf 1.} Omitting the technical details one can rephrase the above
proof in the following way (cf. \cite{W1}).

Let us look first at moves of a real marked point. Consider the
(complex) surface $\cal V$ formed by all the curves $C\in V(\bw')$
and the standard projection $p:{\cal V}\to \Sig$. As it follows
from the proofs in Section \ref{proof}, the branch locus of $p$ is
formed by those curves $C$ which have a cusp, and the branching
index is two. Moreover, when $z_r$ moves in $\R \Sig$ outside the
ramification locus, the curve components $C$ of $p^{-1} (z_r)$,
which correspond to the points of $\Lambda (z_r)\cap\clV(\bw')$,
are subject of a real regular homotopy with an exception of a
finite number of real Morse transformations. Both the regular
homotopies and the Morse transformations preserve the parity of
the number of the solitary nodes (indeed, only one kind, up to
reversing, of regular homotopies is changing this number: two
solitary nodes can come together and then disappear together). At
last, also up to reversing, when $z_r$ crosses the ramification
locus, two nodal curves come together, turn into a curve with a
cusp and then disappear together. On one of them the cusp is
generated by a solitary node, on the other one by a non-solitary
node, and thus they give opposite contributions to the Welschinger
number.

When we move a pair of complex conjugated points, we deal with the
bifurcation locus of the correspondence $C\in
\R\clV(\bw')\cap\Lambda(z_1,z_2)$ between $\R\clV(\bw')$ and $\Sig
\times \Sig$. As it follows from the proofs in
Sections~\ref{proof} and~\ref{proof2}, the trace of the
bifurcation locus on the real part $z_2=\bar z_1$ of $\Sig \times
\Sig$ is a real algebraic variety of real dimension $3$, and
transversal crossings of the bifurcation locus provide the same
bifurcations as in the case of a move of a real point (in fact,
there is even one bifurcation less, since there are no more
bifurcations as in Lemma \ref{lT}). Hence, the Welschinger number
is never changing.

Note that the bifurcations as in Lemma \ref{lT} are a source of
non-invariance of the numbers $W_{r,m}$ in the case of positive
genus, see \cite{IKS}. Note also that, in accordance with the
non-genericity of the complex structure of nodal Del Pezzo
surfaces, the numbers $W_{r,m}$ are not invariant for nodal Del
Pezzo surfaces; the simplest example is splitting out of a
$(-2)$-curve as in \cite{CH}, Proposition 2.6(2b). That is why it
is not clear how to use the Welschinger invariant for enumeration
of curves on nodal Del Pezzo surfaces.

\vskip5pt

\noindent {\bf 2.} One may be interested to extend the results
from $\R$ to other real closed fields. To this end, consider the
$s$-anticanonical models of real unnodal Del Pezzo surfaces of a
given degree~$n$ (and with given $s$ depending on $n$). They form
a Zariski open subset ${\cal H}^0_n(\R)$ in certain irreducible
components of the corresponding Hilbert scheme ${\cal H}_n(\R)$.
Note that ${\cal H}_n(\R)$ is defined over $\Q$, and thus ${\cal
H}^0_n(\R)$ is defined over the field $\R_{\alg}$ of algebraic
real numbers. Consider the family ${\cal F}^0_n(\R)$ of
configuration spaces $\Omega_{r,m}(\Sig)$ over ${\cal H}^0_n(\R)$.
The Welschinger number defines a function~$W$ on ${\cal
F}^0_n(\R)$. This function is semi-algebraic and defined over
$\R_{\alg}$. Thus, $W$ is defined for any real closed field $\K$
and constitutes a semi-algebraic function on the corresponding
family ${\cal F}^0_n(\K)$.

\noindent {\bf Welschinger theorem for real closed fields. } {\it
The function~$W$ is constant on any semi-algebraically connected
component of ${\cal F}^0_n(\K)$.}

{\bf Proof. } For $\K = \R$ the statement follows from \cite{W1},
Theorem 2.1 (for the definition of semi-algebraically connected
components see~\cite{BCR}; in the case $\K = \R$
semi-algebraically connected components coincide with the usual
ones). According to the remarks on~$W$ made above, the
Tarski-Seidenberg principle (see, for example, \cite{BCR}) applies
and gives the statement required.\proofend

Note that the connected components of ${\cal F}^0_n(\R)$ are easy
to describe in topological terms. Namely, $(\Sig, \bw)$ and
$(\Sig', \bw')$, where $\bw \in \Omega_{r,m}(\Sig)$ and $\bw' \in
\Omega_{r,m}(\Sig')$, are in the same component of ${\cal
F}^0_n(\R)$ if $\R\Sigma$ and $\R\Sigma'$ are homeomorphic
(see~\cite{DIK}).

\vskip5pt

\noindent {\bf 3.}~When we were finishing these notes, we have
become acquainted with a new preprint of Welschinger~\cite{W2},
where he proves the invariance of another characteristic, specific
for rational curves in real algebraic convex $3$-folds. The
algebro-geometric approach he develops there (an approach which
uses parametrizations) provides another strategy for a proof of
Theorem~\ref{p1} in the algebraic geometry framework.

{\bf Acknowledgements}. The first and the second authors are
members of Research Training Network RAAG CT-2001-00271; the
second author is as well a member of Research Training Network
EDGE CT-2000-00101; the third author has been supported by the
Bessel research award from the Alexander von Humboldt
Stiftung/Foundation.

{\ncsc Universit\'e Louis Pasteur et IRMA \\[-21pt]

7, rue Ren\'e Descartes, 67084 Strasbourg Cedex, France} \\[-21pt]

{\it E-mail address}: {\ntt itenberg@math.u-strasbg.fr}

\vskip10pt

{\ncsc Universit\'e Louis Pasteur et IRMA \\[-21pt]

7, rue Ren\'e Descartes, 67084 Strasbourg Cedex, France} \\[-21pt]

{\it E-mail address}: {\ntt kharlam@math.u-strasbg.fr}

\vskip10pt

{\ncsc School of Mathematical Sciences \\[-21pt]

Raymond and Beverly Sackler Faculty of Exact Sciences\\[-21pt]

Tel Aviv University \\[-21pt]

Ramat Aviv, 69978 Tel Aviv, Israel} \\[-21pt]

{\it E-mail address}: {\ntt shustin@post.tau.ac.il}

\end{document}